\documentclass[reqno,12pt]{amsart}
\usepackage{amsmath}

\newtheorem{theorem}{Theorem}[section]
\newtheorem{lemma}[theorem]{Lemma}
\newtheorem{corollary}[theorem]{Corollary}
\newtheorem{proposition}[theorem]{Proposition}

\theoremstyle{definition}

\theoremstyle{remark}
\newtheorem{remark}[theorem]{Remark}

\newcommand{\mysection}[1]{\section{#1}
\setcounter{equation}{0}}

\newcommand{\bR}{\mathbb R}

\newcommand{\norm}[1]{\left\Vert#1\right\Vert}
\newcommand{\abs}[1]{\left\vert#1\right\vert}

\renewcommand{\epsilon}{\varepsilon}

\begin{document}
\title[Navier-Stokes equations
]{On the analyticity of solutions to the Navier-Stokes equations
with fractional dissipation}

\author[H. Dong]{Hongjie Dong}
\address[H. Dong]
{School of Mathematics, Institute for Advanced Study,
Einstein Drive, Princeton, NJ 08540, USA}
\email{hjdong@math.ias.edu}

\thanks{This material is based upon work supported by the National Science Foundation under agreement No. DMS-0111298. Any options, findings and conclusions or recommendations expressed in this material are those of the authors and do not necessarily reflect the views of the National Science Foundation.}

\author[D. Li]{Dong Li}
\address[D. Li]
{School of Mathematics, Institute for Advanced Study,
Einstein Drive, Princeton, NJ 08540, USA}
\email{dongli@math.ias.edu}

\subjclass{35Q35}

\keywords{spatial analyticity, Navier-Stokes equations.}

\begin{abstract}
By using a new bilinear estimate,  a pointwise estimate of the
generalized Oseen kernel and an idea of fractional bootstrap, we
show in this note that solutions to the Navier-Stokes equations
 with fractional dissipation are analytic in space variables.
\end{abstract}

\maketitle

\mysection{Introduction}

We are interested in the initial value problem of $d$ dimensional generalized
Navier-Stokes
 equations with fractional dissipation
\begin{align}
                    \label{nse1}
&u_t+u\nabla u+(-\Delta)^{\gamma/2} u+\nabla p=0,\quad \text{div}\, u=0,\\
                    \label{nse2}
&u(0,x)=u_0(x)\quad x\in \bR^d,
\end{align}
where $\gamma\in (1,2]$ is a fixed parameter and the initial data
$u_0$ is in some Banach space to be specified later.

In a well-known paper \cite{kato} Kato proved that for $\gamma=2$
the problem is locally well-posed for $u_0\in L^d$. Kato's method is
based on perturbation theory of the Stokes kernel and is different
from the energy methods used in the
seminal paper \cite{leray}
by Leray. The so called mild solutions are constructed via a fixed
point argument by considering the corresponding integral equations.
Kato's results have been generalized by many authors in various
function spaces. (See, for example, \cite{giga}, \cite{iftimie},
\cite{tataru}, \cite{Planchon}, \cite{Taylor}). With minor
modifications, this method can also be applied to show the local
well-posedness of the generalized Navier-Stokes equation
\eqref{nse1}-\eqref{nse2} with initial data $u_0\in L^{\frac d
{\gamma-1}}$ and the global well-posedness for small data (see
Proposition \ref{prop1}).

In \cite{masuda}, Masuda initiated the study the spatial analyticity
of solutions to the Navier-Stokes equations. The temporal
analyticity was proved by Foias and Temam in an important paper
\cite{foias}. The study of analyticity of the Navier-Stokes
equations was continued by many authors. (See, for example,
\cite{giga2}, \cite{GK}, \cite{LS}, \cite{L00},  \cite{sawada} and
\cite{SM}.) In particular, in a very recent paper \cite{SM}, Miura
and Sawada showed that the solutions by  Koch and Tataru
\cite{tataru} are spatial analytic. A similar smoothness result is
also obtained in a recent preprint \cite{GPS2006} by Germain,
Pavlovi\'c and Staffilani for both the $L^\infty$ and the Carleson
norms.

Usually, spatial analyticity of solutions to the Navier-Stokes
equations is obtained by using either the fixed point argument (see,
for example, \cite{sawada}, \cite{SM} and \cite{GPS2006}) or a
variation of Foias and Temam's method (see, for example,
\cite{L02}). For both methods, one always needs some kind of
smallness assumption on either the initial data or the solution
itself. The main result (Theorem \ref{thm1}) of this paper is that
the spatial analyticity is an intrinsic property of the solutions to
the Navier-Stokes equations. The philosophy is that the existence of
the solution in certain path spaces  implies its analyticity without
any smallness assumption. Our method is based on  new estimates of
the kernels, a semi-group property of the mild solutions and a so
called fractional bootstrap argument.

The remaining part of the article is organized as follows: the main
results are given in the following section. Section \ref{kernelest}
is devoted to a new pointwise estimate of the generalized Oseen kernel and a
corresponding estimate of the generalized heat kernel. In Section
\ref{proofofthm1} we prove Theorem \ref{thm1} mainly by using the
aforementioned fractional bootstrap argument. Finally, the proof of
Corollary \ref{thm2} is given in the last section.

\mysection{The main results}    \label{mainresults}

Define $G(t,x)=G_\gamma(t,x)$ by its Fourier  transform
$\widehat{G_\gamma}(t,\xi)=e^{-t|\xi|^\gamma}$ for $t>0$. Then
$G_\gamma(t,x)$ is the fundamental solution of the linear operator
$\partial_t+(-\Delta)^{\gamma/2}$ and it has the scaling property
\begin{equation}
                \label{eq1.27}
G_\gamma(t,x)=t^{-\frac{d}{\gamma}}G_\gamma(1,t^{-\frac{1}{\gamma}}x).
\end{equation}
It is well-known that \eqref{nse1}-\eqref{nse2} can be rewritten
into an integral equation

\begin{align}
u(t) & =G(t)*u_0-\int_0^t G(t-r,\cdot)*{\mathcal P}(u\nabla
u)(r,\cdot)\,dr,  \nonumber \\
& = G(t)*u_0 - \int_0^t K(t-r,\cdot)*(u\nabla
u)(r,\cdot)\,dr,
                   \label{eq2.33}
\end{align}
where $\mathcal P$ is the Helmholtz projection, and
$K=PG$ is the Oseen kernel (see Section \ref{kernelest} and Proposition
\ref{Oseen}).

For $q\in (\frac{d}{\gamma-1},\infty],T\in(0,\infty]$, introduce the
Banach spaces
$$
X_{q,T}=BC([0,T),L_x^{\frac{d}{\gamma-1}})\cap
\{u\,|\,t^{\alpha}u\in BC((0,T),L_x^q)\},\quad
\alpha=1-\frac{1}{\gamma}-\frac{d}{q\gamma},
$$
with norm
$$
\|u\|_{X_{q,T}}=\max \bigl\{\|u\|_{L_x^\frac{d}{\gamma-1}L_t^\infty([0,T))},
\|t^\alpha u\|_{L_x^q L_t^\infty((0,T))} \bigr\},
$$

The classical Kato's method easily gives the following local well-posedness result.

\begin{proposition}
            \label{prop1}
Assume $u_0$ is in the closure of $\{u\in
C_0^{\infty}(\Omega)\,;\,\text{div}\,u=0\}$ in the scaling invariant
Lebesgue space $L^{\frac{d}{\gamma-1}}(\bR^d)$. Then for any $q\in
(\frac{d}{\gamma-1},\infty]$, \eqref{eq2.33} has a unique solution
in $X_{q,T}$ for some $T\in (0,\infty]$.
\end{proposition}

Here we state our main results of this note.
\begin{theorem}
                    \label{thm1}
Suppose $u$ is a solution to \eqref{eq2.33} in $\bR^d\times (0,T)$
for some $T\in (0,\infty]$ and satisfies  $\|t^\alpha u\|_{L_x^q
L_t^\infty((0,T))}<\infty$ for some $q\in
(\frac{d}{\gamma-1},\infty)$. Then for any $t\in (0,T)$ and $q'\in
[q,\infty]$, we have
\begin{equation}
            \label{eq12.44}
\|D^k u(t,\cdot)\|_{L_x^{q'}}\leq C^{k+1} t^{-\frac{k}{\gamma}-\alpha'}k^k,
\end{equation}
where $\alpha'=1-\frac{1}{\gamma}-\frac{d}{q'\gamma}$ and $C$ is
independent of $k$ and $q'$. Consequently, $u(t,\cdot)$ is spatial
analytic.
\end{theorem}

\begin{remark}
            \label{rem12.42}
If $T=\infty$, estimate \eqref{eq12.44}  implies the decay in time
of higher order Sobolev norms. Furthermore, the radius of
convergence of Taylor's expansion of $u(t,\cdot)$ increases with
time at a rate proportional to $t^{1/\gamma}$.
\end{remark}

\begin{remark}
                    \label{rem1}
In Theorem \ref{thm1} we don't assume any  condition on the initial
data $u_0$. The philosophy here is that the mere existence of the
solution implies its analyticity.
\end{remark}

\begin{remark}
In the case when $\gamma=2$, a similar result is obtained in  a
recent interesting paper \cite{giga2} by Giga and Sawada. Our proof
is more direct, and essentially different from theirs. Moreover, for
general $\gamma\in (1,2]$, their method only gives a less
satisfactory estimate
$$
\|D^k u(t,\cdot)\|_{L_x^q}\leq C^{k+1}
t^{-\frac{k}{\gamma}-\alpha}k^{\frac {2k} {\gamma}},
$$
which doesn't imply the spatial analyticity of $u$ if $\gamma<2$.
\end{remark}

The next corollary is a simple consequence of Proposition
\ref{prop1} and  Theorem \ref{thm1}.
\begin{corollary}
                    \label{thm2}
Suppose $u$ is a solution to \eqref{eq2.33} in $\bR^d\times (0,T)$
for some $T\in (0,\infty]$ and satisfies  $u\in
C([0,T),L_x^{d/(\gamma-1)})$. Then there exists a countable subset
$\Omega$ of $(0,T)$ such that $u(t,\cdot)$ is spatial analytic for
any $t\in (0,T)\setminus \Omega$.
\end{corollary}

%

\mysection{Pointwise estimates of the generalized Oseen kernel}
                                                    \label{kernelest}

\textbf{The generalized Oseen kernel} We will need the following pointwise
estimates of higher derivatives of the generalized Oseen kernel with an explicit
control of constants. The proof of a similar result but with {\em no
control of constants} can be found in \cite{L02}, Proposition 11.1.
We will not use Proposition \ref{Oseen} in its full generality.
However, the estimate itself is of independent interest and we are
not able to find it in the literature.

\begin{proposition} \label{Oseen}
Assume $d\ge 3$ and {$\gamma\in (0,\infty)$}. For $1
\leq j,m < d$ and $t
>0$, the operator $O_{j,m,t} = \frac{1}{\Delta} \partial_j
\partial_m e^{- t \Lambda^{\gamma} }$ is a convolution operator
whose kernel is given by
$$K_{j,m,t} = \frac{1}{t^{d/\gamma }} K_{j,m}(\frac{x}{t^{1/\gamma} }),$$ where
$K_{j,m}$ is a smooth function. There exists  a constant
$C=C(d,\alpha,\gamma)$ such that, for any integer $k\geq 0$ and $ -1
< \alpha \le 1$,
$$
\abs{ (1 + |x|)^{d + k+\alpha } \partial^{k}  \Lambda^{\alpha}
K_{j,m} }  \le C^{k+1} k^k
\mbox{ for all } x\in \bR^d.
$$
\end{proposition}
\begin{proof}
Let $G_{\gamma}(t,y)$ be the generalized heat kernel which satisfies the
scaling property \eqref{eq1.27}.
Consider first the case $0\le  \abs x \le 1$. We have
\begin{align*}
&\max_{ 0 \le \abs x \le 1} \abs{ D^{k+2} \Lambda^{\alpha-2} G_{\gamma}(1,y) } \\
&\le \norm{ D^{k+2} \Lambda^{\alpha-2} G_{\gamma}(1,y)
}_{L^\infty_x} \le \int_{R^d} \abs{\xi}^{k+\alpha} e^{-
\abs{\xi}^\gamma} d\xi \le C^{k+1} k^k.
\end{align*}
Now it is enough to show that $\forall$ $\abs x \ge 1$, we have
$$
D_x^{k+2} \Lambda^\alpha \int \frac 1 { \abs{x-y}^{d-2} } G_{\gamma}(1,y) dy
\le  \frac {C^{k+1} k^k} { \abs{x}^{k+\alpha+d}}.
$$
To this end, write $ x = t^{-1/\gamma} \hat n$, $\hat n \in S^d$,
$y= t^{-1/\gamma} z$, $0<t\le 1$, and we have
\begin{align*}
& \abs{x}^{k+\alpha+d} D_x^{k+2} \Lambda^\alpha \int
\frac 1 { \abs{x-y}^{d-2}} G_\gamma(1,y) dy \\
&= C_1 \abs{x}^{k+\alpha+d} \int \frac 1 {\abs{x-y}^{d-2+\alpha}}
D_y^{k+2} G_\gamma(1,y) dy \\
&=  C_1 \int \frac {1} { \abs{ \hat n - z}^{d-2+\alpha}} D_z^{k+2}
G_\gamma(t, z) dz.
\end{align*}
where $C_1=C_1(d,\gamma,\alpha)$ is another constant.

Now note that as $t\rightarrow 0$, $G_\gamma(t,x) \rightarrow
\delta_0(x)$ where $\delta_0$ is the Dirac distribution on $R^d$.
Therefore it is easy to see that the right-hand side of the above
converges to
$$
\Bigl. D_z^{k+2} ( \frac {1} {\abs {\hat n -z } ^{d-2+\alpha}}) \Bigr|_{z=0}.
$$
Clearly this is bounded by $C^{k+1}\cdot k^k $ for some constant
$C>0$. Remark that this heuristic argument suggests why the optimal
bound on the constants is of the form $C^{k+1} \cdot k^k$.

To complete our argument, we write
\begin{align*}
& \int \frac 1 { \abs{ \hat n -z}^{d-2+\alpha}} D_z^{k+2} G_\gamma(t,z) dz \\
&= \int_{ \abs{ \hat n -z} \le 1/2} \frac 1 { \abs{ \hat n
-z}^{d-2+\alpha}} D_z^{k+2} G_\gamma(t,z) dz\\
&\,\,+\int_{ \abs{ \hat n -z} > 1/2} \frac 1 { \abs{ \hat n -z}^{d-2+\alpha}} D_z^{k+2} G_\gamma(t,z) dz \\
&= \text{I} + \text{II}.
\end{align*}
To estimate I, we use the representation of $G_\gamma(t,z)$ through
the heat kernel \cite{Res95}:
\begin{align}    \label{eq.10.9.05}
G_\gamma(t,z) = \int_0^\infty \frac 1 { (\pi^{\frac 12} s^{\frac 12} t^{\frac 1{\gamma}})^d}
\exp\left\{ - \abs{ \dfrac z  {s^{\frac 12} t^{\frac 1 {\gamma}}} }^2\right\} dF(s),
\end{align}
where $dF(\cdot)$ is a probability measure. This gives us
\begin{align*}
& D_z^{k+2} G_\gamma(t,z) \\
& = \int_0^\infty \frac 1 { (\pi^{\frac 12} s^{\frac 12} t^{\frac
1{\gamma}})^d}
D_z^{k+2} \exp\left\{ - \abs{ \dfrac z  {s^{\frac 12}
t^{\frac 1 {\gamma}}} }^2\right\}\, dF(s) \\
&=  \int_0^\infty \pi^{-\frac d2}\left( t^{\frac 1{\gamma}} s^{\frac
12} \right)^{-(k+d+2)} \cdot (-1)^k \cdot He_{k+2} \left( \frac z
{s^{\frac 12} t^{\frac 1{\gamma}} } \right)\\
&\quad \cdot \exp\left\{ - \abs{ \dfrac z  {s^{\frac 12} t^{\frac 1
{\gamma}}} }^2\right\}\, dF(s).
\end{align*}
where $He_{k+2}(z)$ is the $d$-dimensional Hermite polynomial of
degree $k+2$. We now use the following pointwise estimate of Hermite
polynomials \cite{Ind61}:
\begin{align}   \label{eq.10.10.09}
\abs{ He_n(x) } \le (2^n n!)^{\frac 12} e^{\frac {x^2} 2}.
\end{align}
We then have
\begin{align*}
& \abs{ D_z^{k+2} G_\gamma(t,z) } \\
&\le \int_0^\infty (t^{\frac 1{\gamma}} s^{\frac 12}
)^{-(k+d+2)}\cdot 2^{\frac k 2+1} \cdot ((k+2)!)^{\frac 12} \cdot
\exp\left\{ - \frac 1 2\abs{ \dfrac z  {s^{\frac 12} t^{\frac 1
{\gamma}}} }^2\right\} dF(s).
\end{align*}
Since in case I, $\abs{\hat n -z} \le \frac 12$ and therefore $\abs{z} \ge \frac 12$, we have
\begin{align*}
& \abs{ D_z^{k+2} G_\gamma(t,z) } \\
&\le {C^{k+1}} \cdot ((k+2)!)^{\frac 12}
\displaystyle\int_0^\infty (t^{\frac 1{\gamma}} s^{\frac 12}
)^{-(k+d+2)}
\exp\left\{- \frac 1 {8 \abs{ s^{\frac 1 2} t^{\frac 1{\gamma}}}^2}  \right\} dF(s) \\
&\le  {C^{k+1}} k^k.
\end{align*}
where the last inequality follows from the fact that $dF(s)$ is a probability measure and
the elementary inequality
\begin{equation} \label{eq.Dec.9.57}
\sup_{x>0} x^{k+d+2} e^{-\frac {x^2} 8} \le {C^{k+1}}
\cdot k^{\frac k 2}.
\end{equation}
The estimate of II is similar. By integration by parts and
\eqref{eq.10.9.05} we have,
\begin{align*}
\abs{ \text{II} } &\le \sum_{j=0}^{k+1}  \int_0^\infty (s^{\frac 12} t^{\frac 1{\gamma}})^{-d}
\int_{ \abs{ z-\hat n} = \frac 12}
\abs{ D_z^j ( \frac 1 { \abs{ z-\hat n}^{d-2+\alpha}} ) } \cdot  \\
& \cdot \abs{
 D_z^{k-j+1}
\exp\left\{ - \abs{ \frac z { s^{\frac 12} t^{\frac 1 {\gamma}}} }^2 \right\} }
d\sigma(z) dF(s) +\\
&+ \int_0^\infty (s^{\frac 12} t^{\frac 1{\gamma}} )^{-d}
\int_{\abs{ z-\hat n} > \frac 12} \abs{ D^{k+2}_z ( \frac 1 {\abs{ z-\hat n}^{d-2+\alpha}}) }
\exp\left\{ - \abs{ \frac z { s^{\frac 12} t^{\frac 1 {\gamma}}} }^2 \right\} dz dF(s).
\end{align*}
On $ \abs{ z-\hat n } = \frac 12$, we have $\abs z \ge \frac 1 2$
and therefore by \eqref{eq.10.10.09}
\begin{align*}
& \abs{ D_z^{k-j+1} \exp\left\{ - \abs { \frac z { s^{\frac 12} t^{\frac 1{\gamma}}}}^2 \right\} }\\
&=  (s^{\frac 12} t^{\frac 1{\gamma}} )^{ - (k-j+1)} \cdot \abs{
He_{k-j+1} \left( \frac {z} {s^{\frac 12} t^{\frac 1 {\gamma}}}
\right) }
\cdot \exp\left\{ - \abs { \frac z { s^{\frac 12} t^{\frac 1{\gamma}}}}^2 \right\} \\
&\le  (s^{\frac 12} t^{\frac 1{\gamma}} )^{-(k-j+1)} \cdot 2^{\frac
{k-j+1} 2} ( (k-j+1)! )^{\frac 1 2} \cdot \exp \left\{ - \frac 1 {8
\abs{ s^{\frac 1 2} t^{\frac 1{\gamma}}}^2} \right\}.
\end{align*}
Note also that for $\abs{z-\hat n}>\frac 12$, we have $ \abs{
D^{k+2}_z ( \frac 1 {\abs{ z-\hat n}^{d-2+\alpha}}) } \le {C^{k+1}}
k!$. These estimates together with the elementary inequality
\eqref{eq.Dec.9.57} immediately give us
\begin{align*}
\abs{\text{II} } \le {\sum_{j=0}^{k+1} C^{k+1} j!
(k-j+1)! +C^{k+1}k^k\le C^{k+1} k^k}.
\end{align*}
The proposition is now proved.
\end{proof}

We also need the following lemma.

\begin{lemma} \label{lem.main}
Let $\gamma\in(0,\infty)$, $p\in [1, \infty]$.
Let $k\ge 0$ be an integer and $\epsilon\in (0,1]$.
Then for some constant $C=C(\gamma,\epsilon)>0$,
we have
\begin{equation}
\norm{ D^k\Lambda^{\alpha} G(t,\cdot) }_{L^p_x} \le C^{k+1} k^{\frac
k {\gamma}} t^{- \frac {k+\alpha} {\gamma} - \frac d{\gamma}
(1-\frac 1p)},
\end{equation}
for any $\alpha\in[\epsilon-1,1]$ such that
$k+\alpha\ge \epsilon$ or $k=\alpha=0$. Note that the constant
$C$ can be taken to be independent of $p$.
\end{lemma}
\begin{proof}
This follows from a similar pointwise estimate as in Proposition
\ref{Oseen}. We omit the details.
\end{proof}


\mysection{Proof of Theorem \ref{thm1}}
                \label{proofofthm1}
Since $u$ is divergence
free, we have $u\nabla u=\nabla\cdot(u\otimes u)$. Therefore, by
using integration by parts, the integral equation \eqref{eq2.33} is
equivalent to
\begin{equation}
                    \label{eq2.33b}
u(t)=G(t)*u_0-B(u,u),
\end{equation}
where
$$
B(u,u):=\int_0^t \nabla K(t-s,\cdot)*(u\otimes u)(s,\cdot)\,ds
$$
is a bilinear term. The following lemma is probably known. We
provide a sketched proof for the sake of completeness.
\begin{lemma}
                    \label{lem2.19}
Under the assumptions of Theorem \ref{thm1}, for any $t\in (0,T)$
and $q'\in [q,\infty]$, we have
\begin{equation}
            \label{eq12.44.21}
\|u(t,\cdot)\|_{L_x^{q'}}\leq Ct^{-\alpha'},
\end{equation}
where $\alpha'=1-\frac{1}{\gamma}-\frac{d}{q'\gamma}$ and the
constant $C$ is independent of $q'$.
\end{lemma}
\begin{proof}
We shall use a bootstrap argument. Assume for some positive constant
$C_0$,
$$
\norm{t^\alpha\|u(t,\cdot)\|_{L_x^{q}}}_{L^\infty(0,T)}\leq C_0.$$
We fix a $t\in (0,T)$ and choose $s\in (t/3,2t/3)$ such that
$$
s^\alpha\|u(s,\cdot)\|_{L_x^{q}}\leq C_0.
$$
From \eqref{eq2.33b} and the semigroup property of $G$, it holds that
$$
u(t,\cdot)=G(t-s)*u(s,\cdot)-\int_s^t\nabla K(t-r,\cdot)*(u\otimes u)(r,\cdot)\,dr.
$$
Taking the $L_x^{q'}$ norm on both sides and using Minkowski's
inequality, Young's inequality, H\"older's inequality, Proposition
\ref{Oseen} and Lemma \ref{lem.main}, we get
\begin{align*}
&\|u(t,\cdot)\|_{L^{q'}}\\
&\leq\|G(t-s)*u(s,\cdot)\|_{L^{q'}}+\int_s^t\|\nabla K(t-r,\cdot)*(u\otimes u)(r,\cdot)\|_{L^{q'}}\,dr\\
&\leq \|G(t-s)\|_{L^{r_1}}\|u(s,\cdot)\|_{L^{q}}+\int_s^t\|\nabla K(t-r,\cdot)\|_{L^{r_2}}\|u(r,\cdot)\|^2_{L^{q}}\,dr\\
&\leq C(t-s)^{-\frac d \gamma (1-\frac{1}{r_1})}s^{-\alpha}
+C\int_s^t (t-r)^{-\frac{1}{\gamma}-\frac{d}{\gamma}(1-\frac{1}{r_2})}r^{-2\alpha}\,dr,
\end{align*}
where $r_1$ and $r_2$ satisfy
$$
1+\frac{1}{q'}=\frac{1}{r_1}+\frac{1}{q},\quad
1+\frac{1}{q'}=\frac{1}{r_2}+\frac{2}{q}.
$$
Because $s\in (t/3,2t/3)$,
$$
\|u(t,\cdot)\|_{L^{q'}}\leq Ct^{-\alpha'}+Ct^{-\alpha'}\int_{1/3}^1(1-r)^{-\frac{1}{\gamma}-\frac{d}{\gamma}(1-\frac{1}{r_2})}r^{-2\alpha}\,dr.
$$
Since $q>\frac{d}{\gamma-1}$, the last integral is finite
if
\begin{equation}
            \label{eq3.48}
0\leq \frac{1}{q}-\frac{1}{q'}\leq \frac{1}{2d}(\gamma-1-\frac{d}{q}).
\end{equation}
Hence, \eqref{eq12.44.21} is proved for $q'$ satisfying
\eqref{eq3.48}. A finite iteration of this argument gives
\eqref{eq3.48} for any $q'\in [q,\infty]$. The lemma is proved.
\end{proof}

Now we are ready to prove Theorem \ref{thm1}.\\
\textbf{Proof of Theorem \ref{thm1}}. In this proof we shall denote
by $C_1$ constants which may vary from line to line but does not
depend on $k$ or { $q'$}. Let $N\ge 1$ be an integer
sufficiently large such that
\begin{align*}
1+\frac 1 N + \frac {2d} {Nq}+ \frac d q < \gamma.
\end{align*}
Define a finite sequence of numbers $q_n\ge q$ such that
\begin{align*}
\frac 1 {q_n} = \frac 1 q - \frac n {Nq}, \quad 0\le n \le N.
\end{align*}
Also define for $l=0,1,\cdots, N-1$, $q^\prime\in [q,\infty]$,
\begin{align*}
A(k,l,q') = \norm{ t^{\alpha^\prime+\frac {k+l/N} {\gamma} }
D_x^k \Lambda^{\frac lN} u }_{L_x^{q'} L_t^\infty(0,T)} ,
\end{align*}
and
\begin{align*}
A(k,l) = \sup_{q \le q^\prime \le \infty} A(k,l,q').
\end{align*}
We shall derive a set of recurrent inequalities for $A(k,l)$. To this
end, by using the semigroup proper of $G$, write
$$
u(t, \cdot) = G(\frac t {k+2}) * u(\frac {k+1} {k+2} t)
- \int_{\frac {k+1} {k+2} t}^ t
\nabla K(t-r,\cdot) * ( u\otimes u) (r,\cdot)\,dr.
$$
Call the first term in the above sum linear term and the other
bilinear term. We have four cases. \\
\noindent \texttt{Case 1}: estimate of the linear term, $k\ge 0$,
$1\le l \le N-1$. By Lemma \ref{lem.main} we have
\begin{align*}
& \norm{ t^{\alpha'+\frac {k+l/N} {\gamma}} D_x^k \Lambda^{\frac lN}
\left(
G(\frac t {k+2} ) * u( \frac {k+1} {k+2} t)
\right)
}_{L^{q'}_x L_t^\infty} \\
 &\le  \norm{t^{\frac 1{N\gamma}}
    \Lambda^{\frac 1N} G(\frac t {k+2})}_{L^1_x L^\infty_t}
   \norm{ t^{\alpha'+\frac {k+\frac {l-1} N} {\gamma}} D_x^k
     \Lambda^{\frac {l-1} N} u ( \frac {k+1} {k+2} t)
       }_{L^{q'}_x L^\infty_t}  \\
&\le  C_{\gamma} (k+2)^{\frac 1 {N\gamma}} \cdot
  \left( \frac {k+2} {k+1} \right)^{\alpha' +
   \frac {k+\frac {l-1} N} {\gamma}} A(k,l-1) \\
&\le  C_1 (k+1)^{\frac 1 {N\gamma}} A(k,l-1).
\end{align*}
\noindent
\texttt{Case 2:} estimate of the linear term, $l=0, k\ge 1$. This case
is similar to Case 1 and we have
\begin{align*}
& \norm{ t^{\alpha'+ \frac k{\gamma}} D_x^k \left( G(\frac t{k+2}) *
 u ( \frac {k+1} {k+2} t)  \right)}_{L_x^{q'} L_t^\infty}\\
&=  \norm{ t^{\alpha'+\frac k{\gamma}} D_x \Lambda^{1/N} D_x
\Lambda^{-1} G(\frac t {k+2}) * D_x^{k-1} \Lambda^{\frac {N-1} N} u(
\frac {k+1} {k+2} t)
}_{L_x^{q'} L_t^\infty} \\
&\le  C_1 \cdot (k+1)^{\frac 1 {N\gamma}} A(k-1,N-1).
\end{align*}
\noindent \texttt{Case 3:} estimate of the nonlinear term for $k\ge
0$, $1\le l \le N-1$. Consider $q'\in[q,\infty]$, obviously $\frac 1
{q'} \in [\frac 1 {q_{n+1}}, \frac 1 {q_{n}}]$ for some $0\le n\le
N-1$. For any two functions $f$, $g$, and $0\leq \epsilon<1$,
$2<p<\infty$, the following fractional Leibniz inequality is well
known:
\begin{align*}
\norm{ \Lambda^\epsilon (fg ) }_{L_x^{p/2}}
\le C_p ( \norm{\Lambda^\epsilon f}_{L_x^p} \norm{g}_{L_x^p}
+ \norm{\Lambda^\epsilon g}_{L_x^p} \norm{f}_{L_x^p} ),
\end{align*}
where the constant $C_p$ depends on $p$. In what follows, we shall
only apply the fractional Leibniz inequality when $p=q_n$, $0\le n
\le N-1$. In this way the constants will not depend on
$q'$. Now by Lemma \ref{lem.main}, Proposition \ref{Oseen}, Young's
inequality and fractional Leibniz inequality, we have
\begin{align*}
 &\norm{ t^{\alpha'+ \frac {k+ \frac lN} {\gamma}}
 D_x^k \Lambda^{\frac lN}
 \int_{\frac {k+1} {k+2}t}^t \nabla K(t-r,\cdot) *
 ( u\otimes u) (r)\, dr }_{L_x^{q'}L_t^\infty} \\
&=   \norm{ t^{\alpha'+ \frac {k+ \frac lN} {\gamma}}
 \int_{\frac {k+1} {k+2}t}^t \Lambda^{\frac 1N} \nabla K(t-r,\cdot) *
 \left( D_x^k \Lambda^{\frac {l-1} N} ( u\otimes u) (r)
 \right)\, dr }_{L_x^{q'}L_t^\infty} \\
&\le  C_1 \left( \int_{1- \frac 1{k+2}}^1 (1-r)^{- \frac {1+\frac
1N}{\gamma} - \frac{d}{\gamma} (\frac 2{q_n} - \frac 1 {q^\prime})}
\, dr\right)
\sum_{j=0}^k {\binom k j} \cdot \\
&\quad \cdot \left(1- \frac 1 {k+2} \right)^{-\frac
{k+\frac{l-1}N}{\gamma} -(2-\frac 2{\gamma}- \frac {2d} {\gamma q_n} )}
\cdot
A(j,l-1,q_n) A(k-j,0,q_n) \\
& \le C_1 \sum_{j=0}^k \binom k j
A(j,l-1) A(k-j,0).
\end{align*}
The last integral converges since we have
\begin{align*}
\frac {1+\frac 1 N} {\gamma} + \frac d{\gamma} ( \frac 2 {q_n} -
 \frac 1{q'} )
\le \frac {1+\frac 1 N} {\gamma} + \frac d {\gamma q}+
 \frac d {Nq\gamma} <1.
\end{align*}
\noindent
\texttt{Case 4:} estimate of the nonlinear term for $k\ge 1$ and $l=0$.
This case is similar to Case 3 but slightly trickier. The trick is to write
\begin{align*}
& D_x^k
 \int_{\frac {k+1} {k+2}t}^t \nabla K(t-r,\cdot) *
 ( u\otimes u) (r)\, dr  \\
&=\int_{\frac {k+1} {k+2}t}^t D_x \Lambda^{-\frac {N-1}N}
  \nabla K(t-r,\cdot) *
 D_x^{k-1} \Lambda^{\frac {N-1} N} ( u\otimes u) (r)\,dr.
\end{align*}
Now the rest of the proof follows essentially the same line as in Case 3.
We have
\begin{align*}
&\norm{ t^{\alpha'+ \frac {k} {\gamma}}
 D_x^k
 \int_{\frac {k+1} {k+2}t}^t \nabla K(t-r,\cdot) *
 ( u\otimes u) (r)\, dr }_{L_x^{q'}L_t^\infty} \\
&\le  C_1 \sum_{j=0}^{k-1} {\binom {k-1} j} A(j,N-1) A(k-1-j,0).
\end{align*}
Concluding from the above four cases, and by Lemma \ref{lem2.19}, we
have the following recurrent
inequalities for $A(k,l)$: \\
For $k=0$, $l=0$,
\begin{align*}
A(0,0) \le C;
\end{align*}
For $k\ge 0$, $1\le l\le N$, denote $A(k,N)=A(k+1,0)$, and we have
\begin{align*}
A(k,l) \le C_1 (k+1)^{\frac 1 {N\gamma}} A(k,l-1)
+ C_1 \sum_{j=0}^k {\binom k j} A(j,l-1) A(k-j,0).
\end{align*}
Here $C$ and $C_1$ are constants greater than 1. For $k\ge 0$, $1\le
l\le N$, $0\le j\le k$, denote $n_1=Nj+l-1$, $n_2= N(k-j)$ and
$n=Nk+l$. The following inequality is easy to prove by using
Stirling's formula:
$$
{\binom k j} \le C_1 \left( \frac {n^n} {n_1^{n_1} n_2^{n_2}} \right)^{\frac 1N}.
$$
Then it is not difficult to see that $A(k,l) \le F(Nk+l)$, where
$F(n)$ is a sequence of numbers satisfying
\begin{align*}
F(0) \le C,
\end{align*}
and for $n\ge 1$,
\begin{align*}
F(n) &= C_1 n^{\frac 1 {N\gamma}} F(n-1)\\
&\quad + C_1 \sum_{n_1=0}^{n-1} \frac {n^{n/N} } {n_1^{n_1/N}
(n-1-n_1)^{(n-1-n_1)/N} }
 F(n_1) F(n-1-n_1).
\end{align*}
Clearly $F(n) \le (C_1 C)^{n+1} n^{n/N} G(n)$, where $G(0)=1$ and
$$
G(n) = 2\sum_{n_1=0}^{n-1} G(n_1) G(n-1-n_1).
$$
By the method of formal power series, it is easy to show that for
some constant $C>0$,
$$
G(n) \le C^{n+1}.
$$
This immediately yields that
$$
A(k,l) \le C^{k+1} k^k.
$$
Our theorem is proved.
\mysection{Proof of Corollary \ref{thm2}}
                \label{proofofthm2}
This section is devoted to the proof of Corollary \ref{thm2}. For
any $t\in [0,T)$, by Proposition \ref{prop1}, \eqref{eq2.33} has a
unique solution $\bar u\in X_{q,\epsilon_t}$ with initial data $\bar
u_0=u(t,\cdot)$ for some $q\in (\frac{d}{\gamma-1},\infty]$ and
$\epsilon_t\in (0,T-t)$.  By Theorem \ref{thm1}, $\bar u(s,\cdot)$
is spatial analytic for $s\in (0,\epsilon_t)$. The same
 is true for $u(t+s,\cdot)$ because of the
following uniqueness result of the mild solution to
\eqref{nse1}-\eqref{nse2}.
\begin{lemma}
                \label{lem5.21}
The mild solution to \eqref{eq2.33} in $C([0,T_1), L_x^{d/(\gamma-1)})$ is unique
 for any $T_1\in(0,\infty]$.
\end{lemma}
Indeed, for $\gamma=2$ this lemma is proved in \cite{Monniaux}. The
proof there can be easily modified  to cover the case $\gamma\in
(1,2]$. We leave the details to interested readers.

Denote
$$
\Omega=\{t\in (0,T)\,|\,u(t,\cdot)\, \text{is not spatial analytic}\}.
$$
For any $t\in \Omega$, we choose a rational number $q_t$ in
$(t,t+\epsilon_t)$. It is easy to see that for different such $t_1$
and $t_2$, the corresponding $q_{t_1}$ and $q_{t_2}$ can not
coincide with each other. Therefore, $\Omega$ is at most a countable
set, and the corollary is proved.

\noindent \textbf{Concluding remark}  with some minor modifications
our method also applies to the periodic boundary condition case. We
leave the details to interested readers.


\end{document}